\documentclass[12pt]{article}

\usepackage{latexsym}
\usepackage{calc}
\usepackage{amssymb, amsmath, amsfonts, listings, mathrsfs} 

\newcounter{hours}\newcounter{minutes}

\def\argmin{\mathop{\rm argmin}}

\def\nr{\par \noindent}

\def\Def{\stackrel{\mathrm{def}}{=}}

\def\dom{{\rm dom \,}}
\def\beq{\begin{equation}}
\def\eeq{\end{equation}}

\def\E{\mathbb{E}}

\def\BI{\begin{itemize}}
	\def\EI{\end{itemize}}

\newcommand{\SetEQ}{\setcounter{equation}{0}}
\newcommand{\refLE}[1]{\ensuremath{\stackrel{\eqref{#1}}{\leq}}}
\newcommand{\refEQ}[1]{\ensuremath{\stackrel{\eqref{#1}}{=}}}
\newcommand{\refGE}[1]{\ensuremath{\stackrel{\eqref{#1}}{\geq}}}

\newtheorem{theorem}{Theorem}
\newtheorem{lemma}{Lemma}
\newtheorem{corollary}{Corollary}

\newtheorem{assumption}{Assumption}
\newtheorem{definition}{Definition}

\newtheorem{example}{Example}
\newtheorem{remark}{Remark}
\newcommand{\proof}{\bf Proof: \rm \nr}
\newcommand{\qed}{\hfill $\Box$ \nr \medskip}
\newcommand{\half}{\mbox{${1 \over 2}$}}

\def\ba{\begin{array}}
	\def\ea{\end{array}}
\def\beann{\begin{eqnarray*}}
	\def\eeann{\end{eqnarray*}}
\def\bea{\begin{eqnarray}}
\def\eea{\end{eqnarray}}

\def\BT{\begin{theorem}}
	\def\ET{\end{theorem}}
\def\BL{\begin{lemma}}
	\def\EL{\end{lemma}}
\def\BC{\begin{corollary}}
	\def\EC{\end{corollary}}
\def\BE{\begin{example}}
	\def\EE{\end{example}}
\def\BD{\begin{definition}}
	\def\ED{\end{definition}}
\def\BR{\begin{remark}}
	\def\ER{\end{remark}}
\def\BAS{\begin{assumption}}
	\def\EAS{\end{assumption}}
\def\BI{\begin{itemize}}
	\def\EI{\end{itemize}}

\def\BMP{\begin{minipage}{9.5cm}}
	\def\EMP{\end{minipage}}
\def\MPT{\begin{minipage}{11.5cm}}
	\def\EPT{\end{minipage}}

\def\la{\langle}
\def\ra{\rangle}

\def\QF{\hspace{5ex} \Box}
\def\QR{\hfill \Box}

\title{
	\textbf{Gradient Regularization of Newton Method with Bregman Distances} \thanks{This project has received funding from the European Research Council
		(ERC) under the European Union's Horizon 2020 research and innovation
		programme (grant agreement No. 788368).} }

\author{Nikita Doikov \thanks{Institute of Information and Communication Technologies,
		Electronics and Applied Math. (ICTEAM), Catholic University of Louvain (UCL). E-mail:
		Nikita.Doikov@uclouvain.be. ORCID: 0000-0003-1141-1625.} 
	\and Yurii Nesterov
	\thanks{Center for Operations Research and Econometrics (CORE),
		Catholic University of Louvain (UCL). E-mail:
		Yurii.Neterov@uclouvain.be. ORCID: 0000-0002-0542-8757.}}

\date{ December 6, 2021 }

\begin{document}
	
	\maketitle
	
	\abstract{
		In this paper, we propose a first second-order scheme based on arbitrary non-Euclidean norms, incorporated by Bregman distances. They are introduced directly in the Newton iterate with regularization parameter
proportional to the square root of the norm of the current gradient. For the basic scheme, as applied to the composite optimization problem, we establish the global convergence rate of the order $O(k^{-2})$ both in terms of the functional residual and in the norm of subgradients. Our main assumption on the smooth part of the objective is Lipschitz continuity of its Hessian. For uniformly convex functions of degree three, we justify global linear rate, and for strongly convex function we prove the local superlinear rate of convergence. Our approach can be seen as a relaxation of the Cubic Regularization of the Newton method \cite{nesterov2006cubic},
which preserves its convergence properties, while the auxiliary subproblem at each iteration is simpler.
We equip our method with adaptive line search procedure for choosing the regularization parameter. We propose also an accelerated scheme with  convergence rate $O(k^{-3})$, where $k$ is the iteration counter. 
	}
	
	\vspace{10ex}\noindent
	
	\thispagestyle{empty}
	
	\newpage\setcounter{page}{1}

	\section{Introduction}

	The classical Newton's method is a powerful tool for solving various optimization problems 
	and for dealing with ill-conditioning. The practical implementation of this method
	for solving unconstrained minimization problem $\min\limits_{x} f(x)$ can be written as
	follows:
	$$
	\ba{rcl}
	x_{k + 1} & = & x_k - \alpha_k \bigl( \nabla^2 f(x_k) \bigr)^{-1} \nabla f(x_k), \qquad k \geq 0,
	\ea
	$$
	where $0 < \alpha_k \leq 1$ is a damping parameter.
	However, this approach has two serious drawbacks.
	Firstly, the next point is not well-defined when the Hessian is not strictly positive-definite.
	And secondly, while the method has a very fast local quadratic convergence, it is difficult to establish any \textit{global}  properties for this process.
	Indeed, for $\alpha_k = 1$ (the classical pure Newton method), there are known examples of problems
	for which the method does not converge globally (see, e.g., Example 1.4.3 in~\cite{doikov2021new}).
	For the damped Newton method with line search, it is possible to prove some global convergence rates.
	But, typically, they are worse than the rates of the classical Gradient Method~\cite{nesterov2018lectures}.
	
	A breakthrough in the second-order optimization theory was made after~\cite{nesterov2006cubic},
	where the Cubic Regularization of the Newton method was presented together 
	with its global convergence properties. The main standard assumption is that the Hessian 
	of the objective is Lipschitz continuous with some parameter $L_2 \geq 0$:
	$$
	\ba{rcl}
	\| \nabla^2 f(x) - \nabla^2 f(y) \| & \leq & L_2 \|x - y\|, \qquad \forall x, y,
	\ea
	$$
	ensuring the \textit{global upper approximation} of our function
	formed by the second-order Taylor polynomial augmented by the third power of the norm.
	The next point is then defined as the minimum of the upper model:
	\beq \label{met-IntroCNM}
	\ba{rcl}
	x_{k + 1} & = & \argmin\limits_{y} \Bigl[ 
	\la \nabla f(x_k), y - x_k \ra + \frac{1}{2} \la \nabla^2 f(x_k)(y - x_k), y - x_k \ra \\
	\\
	& & \qquad \quad \;\; + \; \frac{L_2}{6}\|y - x_k\|^3
	\Bigr].
	\ea
	\eeq
	Till now, this idea has a full theoretical justification only for the Euclidean norm $\| \cdot \|$. 
	In this case, the solution to the auxiliary minimization problem~(\ref{met-IntroCNM})
	does not have a closed form expression, 
	but it can be found by solving a one-dimensional nonlinear equation 
	and by using the standard factorization tools of Linear Algebra.
	However, even in the Euclidean case, the presence of the cubic term in the objective prevents the usage of gradient-type methods 
	(like the conjugate gradients, etc.). This drawback does not allow the application of method (\ref{met-IntroCNM}) 
	to large-scale problems. 
	
	In this paper, we show how to avoid these restrictions.
	Namely, we will show that it is possible to use a \textit{quadratic regularization}
	of the Taylor polynomial with a properly chosen coefficient that depends only on the current iterate.
	In the simplest form, one iteration of our method is as follows:
	\beq \label{met-IntroQ}
	\ba{rcl}
	x_{k + 1} & = & x_k - \bigl( \nabla^2 f(x_k) + A_k I \bigr)^{-1} \nabla f(x_k),
	\ea
	\eeq
	where
	\beq \label{met-IntroAlpha}
	\ba{rcl}
	A_k & = & \sqrt{\frac{L_2}{3} \| \nabla f(x_k)\|}.
	\ea
	\eeq
	We see that it is very easy for implementation, since it requires only \textit{one}
	matrix inversion, the very standard operation of Linear Algebra.
	At the same time, this subproblem is now suitable for the classical Congugate Gradient method as well.\footnote{$^)$ When this paper was already finished, we discovered that this idea was recently proposed by K. Mishchenko \cite{mishchenko2021regularized} for solving unconstrained minimization problem with smooth objective. As compared to his work, our main advances consist in the usage of Bregman distances, composite form of optimization problem, linear rate of convergence for uniformly convex functions, and developments of accelerated variant of the method.}$^)$
	
	It appears that for the optimization process~\eqref{met-IntroQ},\eqref{met-IntroAlpha},
	we can establish the global convergence guarantees of the same type as for the Cubic Newton method~\eqref{met-IntroCNM}. Namely, we prove the global rate of the order
	$O(1/k^2)$ in terms of the functional residual and in terms of the subgradient norm 
	for the general convex functions.
	This is much faster than the standard $O(1/k)$-rate of the Gradient Method.
	Moreover, for the uniformly convex functions of degree three, we prove the global 
	linear rate. For the strongly convex functions we establish a local superlinear convergence.
	
	\vspace{1ex}\noindent
	{\bf Contents.}
	In this paper, we consider optimization problems in a general composite form. We can work with arbitrary (possibly non-Euclidean) norms using the framework of Bregman distances.
	
	In Section~\ref{sc-GReg}, we present the main properties of one iteration of the scheme.
	We study the convergence properties of the basic process in Section~\ref{sc-MinProc}.
	In Section~\ref{sc-Grad}, we establish convergence rates for the norm of the gradient. A line search procedure for our scheme is discussed in Section~\ref{sc-LS}. In Section~\ref{sc-Acc}, we consider an accelerated method based on the iterations of the basic process and justify its global complexity of the order $\tilde{O}(\epsilon^{-1/3})$ assuming Lipschitz continuity of the Hessian of the smooth part of the objective function.

	\vspace{1ex}\noindent
	{\bf Notation.}  
	Let us fix a finite-dimensional real vector space $\E$.
	Our goal is to solve the following {\em Composite Minimization Problem}
	\beq\label{prob-Uni}
	F^* \; = \; \min\limits_{x \in \dom \psi} \bigl[ F(x) \; \Def \; f(x) + \psi(x) \bigr],
	\eeq
	where $\psi(\cdot)$ is a {\em simple} closed convex function 
	with $\dom \psi \subseteq \E$, and $f(\cdot)$ is a convex and two times continuously differentiable function. 
	
	We measure distances in $\E$ by a general norm $\| \cdot \|$.
	Its dual space is denoted by $\E^{*}$. It is a space of all linear functions on $\E$, for which we define the norm in the standard way:
	$$
	\ba{rcl}
	\| g \|_* & = & \max\limits_{x \in \E} \{ \; \la g, x \ra: \; \| x \| \leq 1 \; \}, 
	\qquad g \in \E^*.
	\ea
	$$
	Using this norm, we can define an induced norm for a self-adjoint linear operator $B: \E \to \E^*$ as follows:
	$$
	\ba{rcl}
	\| B \| & = & \max\limits_{x \in \E} \{ | \la B x, x \ra | : \; \| x \| \leq 1 \}.
	\ea
	$$
	We can also define the bounds of its spectrum as the best values $\lambda_{\min}(B)$ and $\lambda_{\max}(B)$ satisfying conditions
	$$
	\ba{rcl}
	\lambda_{\min}(B) \| x \|^2 & \leq & \la B x, x \ra \; \leq \; \lambda_{\max}(B) \| x \|^2, \quad \forall x \in \E.
	\ea
	$$
	
	Our optimization schemes will be based on some scaling function $d(\cdot)$, which we assume to be a strongly convex function with Lipschitz-continuous gradients:
	\beq\label{def-DStrong}
	\ba{rcl}
	d(y) & \geq & d(x) + \la \nabla d(x), y - x \ra + {\sigma \over 2} \| y - x \|^2,
	\ea
	\eeq
	\beq\label{def-DLip}
	\ba{rcl}
	\| \nabla d(x) - \nabla d(y) \|_* & \leq & \| x - y \|,
	\ea
	\eeq
	where $\sigma \in (0,1]$ and the points $x, y \in \dom \psi$ are arbitrary. For twice-differentiable scaling functions, this condition can be characterized by the following bounds on the Hessian:
	$$
	\ba{rcl}
	\sigma \|h\|^2 & \leq & \la \nabla^2 d(x)h, h \ra
	\;\; \leq \;\; \|h\|^2, \qquad \forall x \in \dom \psi, \; h \in \E.
	\ea
	$$
Using this function, we define the following {\em Bregman distance}:
	\beq\label{def-Breg}
	\ba{rcl}
	\rho(x,y) & = & \beta_d(x,y) \;\; \Def \;\; 
	d(y) - d(x) - \la \nabla d(x), y - x \ra, \quad x, y \in \dom \psi.
	\ea
	\eeq
	
	The standard condition for the smooth part of the objective function in problem (\ref{prob-Uni}) is Lipschitz continuity of the Hessians:
	\beq\label{eq-Lip}
	\| \nabla^2 f(x) - \nabla^2 f(y) \| \;\; \leq \;\; L_2 \| x - y \|, \qquad \forall x, y \in \dom \psi.
	\eeq
	This inequality has the following consequences, which are valid for all $x,y \in \dom \psi$:
	\beq\label{eq-GLip}
	\| \nabla f(y) - \nabla f(x) - \nabla^2 f(x)(y-x) \|_* \;\; \leq \;\; \half L_2 \| y - x \|^2, 
	\eeq
	\beq\label{eq-FLip}
	\ba{c}
	|f(y) - f(x) - \la \nabla f(x), y - x \ra + \half \la \nabla^2 f(x)(y-x), y-x \ra | 
	\;\; \leq \;\; \frac{1}{6} L_2 \| y - x \|^3.
	\ea
	\eeq
	
	\section{Gradient regularization}\label{sc-GReg}
	\SetEQ

	Our main iteration at some point $\bar x \in \dom \psi$ with a step-size $A>0$ is defined as follows: 
	\beq\label{def-TA}
	\ba{rl}
	T_A(\bar x) \Def \arg\min\limits_{y \in \dom \psi} & \Big[ \quad M_A(\bar x, y) \Def f(\bar x) + \la \nabla f(\bar x), y - \bar x \ra \\
	\\
	& + \half \la \nabla^2 f(\bar x)(y - \bar x), y - \bar x \ra + A \rho(\bar x, y) + \psi(y) \quad \Big].
	\ea
	\eeq
	The solution to this problem $T = T_A(\bar x)$ is characterized by the following variational principle:
	\beq\label{eq-TVar}
	\ba{c}
	\la \nabla f(\bar x) + \nabla^2f(\bar x)(T-\bar x) + A (\nabla d(T) - \nabla d(\bar x)), y - T \ra\\
	\\
	+ \psi(y)\; \geq \; \psi(T), \quad y \in \dom \psi.
	\ea
	\eeq
	Thus, defining $\psi'(T) = - \nabla f(\bar x) - \nabla^2f(\bar x)(T-\bar x) - A (\nabla d(T) - \nabla d(\bar x))$, we see that $\psi'(T) \in \partial \psi(T)$. Consequently,
	\beq\label{def-SubT}
	\ba{c}
	F'(T) \; = \; \nabla f(T) + \psi'(T)\\
	\\
	= \; \nabla f(T) - \nabla f(\bar x) - \nabla^2f(\bar x)(T-\bar x) - A (\nabla d(T) - \nabla d(\bar x)) \; \in \; \partial F(T).
	\ea
	\eeq
Note that this is a very special way of selecting subgradient of a possibly nonsmooth function $F(\cdot)$, which allows $\| F'(T) \|_*$ approach zero.
	
	Denote $M_A(\bar x) = M_A(\bar x, T_A(\bar x)) \leq M_A(\bar x, \bar x) = F(\bar x)$. Let us prove the following fact.
	\BL\label{lm-TSize}
	For all $y \in \dom \psi$ and $T = T_A(\bar x)$, we have
	\beq\label{eq-MLow}
	\ba{rcl}
	M_A(\bar x, y) & \geq & M_A(\bar x) + \half \la \nabla^2 f(\bar x) (y - T), y - T \ra + \half \sigma A\| y - T \|^2.
	\ea
	\eeq
	Moreover,
	\beq\label{eq-TSize}
	\ba{rcl}
	\| T_A(\bar x) - \bar x \| & \leq & {1 \over \sigma A} \| F'(\bar x) \|_*,
	\ea
	\eeq
	where $F'(\bar x) = \nabla f(\bar x) + \psi'(\bar x)$ and $\psi'(\bar x)$ is an arbitrary element of $\partial \psi(\bar x)$. 
	\EL
	\proof
	For optimization problem in (\ref{def-TA}), define the scaling function
	$$
	\ba{rcl}
	\xi(x) & = & \half \la \nabla^2f(\bar x) x, x \ra + A d(x).
	\ea
	$$
	Note that the objective function in this problem is strongly convex relatively to $\xi(\cdot)$ with constant one. Therefore,
	$$
	\ba{rcl}
	M_A(\bar x, \bar x) - M_A(\bar x) & \geq & \beta_{\xi}(T,y) \; = \; \half \la \nabla^2 f(\bar x) (y - T), y - T \ra + A \beta_d(T,y)\\
	\\
	& \refGE{def-DStrong} & \half \la \nabla^2 f(\bar x) (y - T), y - T \ra + \half \sigma A \| y - T \|^2.
	\ea
	$$
	In order to prove (\ref{eq-TSize}), note that
	$$
	\ba{rcl}
	M_A(\bar x) & \geq & F(\bar x) + \min\limits_{y \in \dom \psi} \Big[ \la F'(\bar x), y - \bar x \ra + \half \sigma A \| y - \bar x \|^2 \Big]\\
	\\
	& \geq & F(\bar x) + \min\limits_{y \in \E} \Big[ \la F'(\bar x), y - \bar x \ra + \half \sigma A \| y - \bar x \|^2 \Big]\\
	\\
	& = & F(\bar x) - {1 \over 2 \sigma A} \| F'(\bar x) \|_*^2.
	\ea
	$$
	Since $M_A(\bar x, \bar x) = F(\bar x)$, we get (\ref{eq-TSize}) from (\ref{eq-MLow}) with $y = \bar x$.
	\qed
	
	In what follows, the parameter $A$ in the optimization problem (\ref{def-TA}) is chosen as
	\beq\label{def-AH}
	\ba{rcl}
	A \; = \; A_H(\bar x) \; = \; {1 \over \sigma} \sqrt{{H \over 3 } \| F'(\bar x) \|_*},
	\ea
	\eeq
	where $H > 0$ is an estimate of the Lipschitz constant $L_2$ in (\ref{eq-Lip}). This choice is explained by the following result.
	\BC\label{cor-AH}
	For $A = A_H(\bar x)$, we have 
	\beq\label{eq-ABound}
	\ba{rcl}
	H \| T_A(\bar x) - \bar x \| & \leq & 3 \sigma A.
	\ea
	\eeq
	\EC
	\proof
	Indeed, this is a simple consequence of inequality (\ref{eq-TSize}) and definition (\ref{def-TA}).
	\qed
	
	Let us relate the optimal value of the auxiliary problem (\ref{def-TA}) with the cubic over-approximation (\ref{eq-FLip}).
	\BL\label{lm-Approx}
	Let $A = A_H(\bar x)$ and $T = T_A(\bar x)$. Assume that for some $H > 0$ the following condition is satisfied:
	\beq\label{eq-FTUp}
	\ba{rcl}
	f(T) & \leq & f(\bar x) + \la \nabla f(\bar x), T - \bar x \ra + \half \la \nabla^2f(\bar x)(T - \bar x), T - \bar x \ra  + {H \over 6} \| T - \bar x \|^3.
	\ea
	\eeq
	Then
	\beq\label{eq-DecF}
	\ba{rcl}
	F(\bar x) - F(T) & \geq & \half \la \nabla^2 f(\bar x) (T - \bar x), T - \bar x \ra + \half \sigma A\| T - \bar x \|^2.
	\ea
	\eeq
	\EL
	\proof
	Indeed,
	$$
	\ba{rcl}
	f(T) & \refLE{eq-FTUp} & M_A(\bar x) - A \rho(\bar x,T) -\psi(T) +{H \over 6} \| T - \bar x \|^3\\
	\\
	& \refLE{def-DStrong} &  M_A(\bar x) -\psi(T) +{H \over 6} \| T - \bar x \|^3 - \half \sigma A \| T - \bar x \|^2\\
	\\
	& \refLE{eq-ABound} &  M_A(\bar x) -\psi(T).
	\ea
	$$
	Thus, $F(T) \leq M_A(\bar x)$ and (\ref{eq-DecF}) follows from (\ref{eq-MLow}) with $y = \bar x$.
	\qed
	
	Finally, we need to estimate the norm of subgradient at the new point.
	\BL\label{lm-GNorm}
	Let $A = A_H(\bar x)$ and $T = T_A(\bar x)$.  Then
	\beq\label{eq-GNorm}
	\ba{rcl}
	\| F'(T) \|_* & \leq &  \sigma A \left( \sigma^{-1} + {3 L_2 \over 2 H} \right) \| T - \bar x \| \; \leq \; c \| F'(\bar x ) \|_*,
	\ea
	\eeq
	where 
	$$
	\ba{rcl}
	c & \Def &  \sigma^{-1} + {3 L_2 \over 2 H} .
	\ea
	$$
	\EL
	\proof
	Indeed,
	$$
	\ba{c}
	\| F'(T) \|_* \refEQ{def-SubT}   \| \nabla f(T) - \nabla f(\bar x) - \nabla^2f(\bar x)(T-\bar x) - A (\nabla d(T) - \nabla d(\bar x)) \|_*\\
	\\
	\refLE{eq-GLip}   \half L_2 \| T - \bar x \|^2 + A \| \nabla d(T) - \nabla d(\bar x) \|_* \; \refLE{def-DLip} \; \half L_2 \| T - \bar x \|^2 + A \| T - \bar x \|_*\\
	\\
	\refLE{eq-ABound}   A \left( 1 + {3 \sigma L_2 \over 2 H} \right) \| T - \bar x \|.
	\ea
	$$
	This is the first inequality in (\ref{eq-GNorm}). For the second one, we can continue as follows:
	$$
	\ba{rcl}
	\| F'(T) \|_* & \refLE{eq-ABound} &  \left( 1 + {3 \sigma L_2 \over 2 H} \right) \cdot {3 \sigma A^2 \over H} \; \refEQ{eq-ABound} \; c  \| F'(\bar x ) \|_*. \QF
	\ea
	$$
	
	Now we can prove the main theorem of this section.
	\BT\label{th-Main}
	Let $A = A_H(\bar x)$ and $T = T_A(\bar x)$.  If for this point relation (\ref{eq-FTUp}) is valid, then
	\beq\label{eq-Prog}
	\ba{rcl}
	F(\bar x) - F(T) & \geq &
	{1
		\over 2 c^2 }  \sqrt{3 \over H} \cdot { \| F'(T) \|_*^2 \over \| F'(\bar x) \|_*^{1/2}}
	.
	\ea
	\eeq
	\ET
	\proof
	We only need to insert in (\ref{eq-DecF}) the first inequality of (\ref{eq-GNorm}) and definition (\ref{def-AH}).
	\qed
	
	\section{Properties of the minimization process}\label{sc-MinProc}
	\SetEQ
	
	Now we can analyze the following minimization process.
	\beq\label{met-GNewt}
	\ba{|l|}
	\hline \\
	\mbox{{\bf Initialization.} Choose $H \geq L_2$, $x_0 \in \dom \psi$, and $F'_0 \in \partial F(x_0)$.}\\
	\\
	\mbox{{\bf $k$th iteration ($k \geq 0$).} 1). Set $g_k = \| F'_k \|_*$ and $A_k = {1 \over \sigma}\sqrt{{H \over 3} g_k}$.}\\
	\\
	\mbox{2). Compute $x_{k+1} = T_{A_k}(x_k)$ and define}\\
	\\
	F'_{k+1} = \nabla f(x_{k+1}) - \nabla f(x_k) - \nabla^2f(x_k)(x_{k+1}-x_k) - A_k(\nabla d(x_{k+1}) - \nabla d(x_k)) .\\
	\\
	\hline
	\ea
	\eeq
	
	Let us introduce the distance to the initial level set:
	$$
	\ba{rcl}
	D & = & \sup\limits_{x \in \dom \psi} \{ \| x - x^* \|: \; F(x) \leq F(x_0) \},
	\ea
	$$
	which we assume to be bounded: $D < +\infty$. We can prove the following convergence rate 
	for method~\eqref{met-GNewt}.
	
	\BT\label{th-Rate}
	Let $H \geq L_2$ and $F(x_k) - F^* \geq \epsilon$ for some $k \geq 0$. Then, 
	\beq\label{eq-Rate}
	\ba{rcl}
	{1 \over [F(x_k) - F^*]^{1/2}} & \geq & {1 \over [F(x_0) - F^*]^{1/2}} + 
	{1 \over 4 c^2} \sqrt{3 \over H D^3} \left( k - \ln  {(F(x_0) - F^*)  \| F'(x_0) \|_*^{1/2}D^{1/2} \over \epsilon^{3/2}} \right).
	\ea
	\eeq
	\ET
	\proof
	Denote $F_k = F(x_k) - F(x^*)$ and $g_k = \| F'(x_k) \|_*$. Thus, $F_k \leq D g_k$.  Note that
	$$
	\ba{rcl}
	{1 \over F_{k+1}^{1/2}} - {1 \over F_k^{1/2}} & = & { F_{k}^{1/2} - F_{k+1}^{1/2} \over F_{k}^{1/2} F_{k+1}^{1/2}} \; = \; {F_k - F_{k+1} \over F_{k}^{1/2} F_{k+1}^{1/2} (F_{k}^{1/2} + F_{k+1}^{1/2})} \; \geq \; {F_k - F_{k+1} \over 2 F_{k} F_{k+1}^{1/2}}.
	\ea
	$$
	Since for all $k \geq 1$, the subgradients of $\psi(\cdot)$ are defined by the rule (\ref{def-SubT}), we can use the results of Section \ref{sc-GReg}.   
	We can continue as follows:
	$$
	\ba{rcl}
	{1 \over F_{k+1}^{1/2}} - {1 \over F_k^{1/2}} & \refGE{eq-Prog} & {\sqrt{3} g_{k+1}^2 \over 4 \sqrt{H} c^2 g_k^{1/2} F_{k} F_{k+1}^{1/2}} \; \geq \; {\sqrt{3}g_{k+1}^{1/2} F_{k+1} \over 4 \sqrt{H} c^2 g_k^{1/2} F_{k} D^{3/2}}\; = \; {g_{k+1}^{1/2} F_{k+1} \over 4c^2 g_k^{1/2}  F_{k}} \sqrt{3 \over HD^3} .
	\ea
	$$
	Summing up these inequalities, we get
	\beq\label{eq-Int}
	\ba{rcl}
	{1 \over F_k^{1/2}} - {1 \over F_0^{1/2}} & \geq & {1 \over 4 c^2} \sqrt{3 \over HD^3} \sum\limits_{i=0}^{k-1} {F_{i+1} g_{i+1}^{1/2} \over F_i g_i^{1/2}} \; \geq \; {k \over 4c^2} \sqrt{3 \over HD^3} \left({F_{k} g_{k}^{1/2} \over F_0 g_0^{1/2}}\right)^{1/k}\\
	\\
	& \geq & {k \over 4c^2} \sqrt{3 \over HD^3} \left({\epsilon^{3/2} \over F_0 g_0^{1/2}D^{1/2}}\right)^{1/k}.
	\ea
	\eeq
	Since
	$$
	\ba{rcl}
	\left({\epsilon^{3/2} \over F_0 g_0^{1/2}D^{1/2}}\right)^{1/k} & = &\exp\Big( - {1 \over k} \ln  {F_0 g_0^{1/2}D^{1/2} \over \epsilon^{3/2}} \Big) \; \geq 1 - {1 \over k} \ln  {F_0 g_0^{1/2}D^{1/2} \over \epsilon^{3/2}},
	\ea
	$$
	we obtain inequality (\ref{eq-Rate}).
	\qed
	\BC
	The second condition of Theorem \ref{th-Rate} can be valid only for
	\beq\label{eq-Comp}
	\ba{rcl}
	k & \leq & 4 c^2 \sqrt{ H D^3 \over 3 \epsilon} + \ln  {(F(x_0) - F^*)  \| F'(x_0) \|_*^{1/2}D^{1/2} \over \epsilon^{3/2}}.
	\ea
	\eeq
	\EC
	\BR\label{rm-H}
	The right-hand side of inequality (\ref{eq-Comp}) can be used for defining the optimal value of parameter $H$. Indeed, it can be chosen as a minimizer of the following function:
	$$
	\ba{c}
	2 \ln(2H\sigma^{-1} + 3L_2) - {3 \over 2} \ln H.
	\ea
	$$
	This gives us
	\beq\label{eq-OptH}
	\ba{rcl}
	H_* & = & {9 \over 2} L_2 \sigma.
	\ea
	\eeq
	In this case,
	\beq\label{eq-OptComp}
	\ba{rcl}
	4c^2\sqrt{ H_* D^3 \over 3 \epsilon} & = & {64 \over 9 \sigma} \sqrt{ 3 L_2 D^3 \over 2 \epsilon \sigma} \; < \; 8.71 \sqrt{ L_2 D^3 \over \epsilon \sigma^3}. \QF
	\ea
	\eeq
	\ER
	
	Let us estimate now the performance of method (\ref{met-GNewt}) on uniformly convex functions. Consider the case when function $F(\cdot)$ is uniformly convex of degree three:
	\beq\label{def-UC3}
	\ba{rcl}
	F(y) & \geq & F(x) + \la F'(x), y - x \ra + {\sigma_3 \over 3} \| y - x \|^3, \quad x, y \in \dom \psi.
	\ea
	\eeq
	For the composite $F(\cdot)$, this property can be ensured either by its smooth component $f(\cdot)$, or by the general component $\psi(\cdot)$. In the latter case, it is not necessary to coordinate this assumption with the smoothness condition (\ref{eq-Lip}).
	
	In our analysis, we need the following straightforward consequence of definition (\ref{def-UC3}):
	\beq\label{eq-UBound}
	\ba{rcl}
	F(x) - F^* & \leq & {2 \over 3 \sqrt{\sigma_3}} \| F'(x) \|_*^{3/2}, \quad x \in \dom \psi.
	\ea
	\eeq
	\BT\label{th-UC3}
	Let $F(\cdot)$ satisfies condition (\ref{def-UC3}). Then for all $k \geq 0$ we have
	\beq\label{eq-RateUC3}
	\ba{rcl}
	F(x_k) - F^* & \leq & D\| F'(x_0) \|_* \cdot \exp\left( - {k \ln(1+S) \over c^{1/2} + \half \ln(1+S)} \right),
	\ea
	\eeq
	where $S = {3 \sqrt{3}  \over 4 c^{3/2} } \sqrt{\sigma_3 \over H}$.
	\ET
	\proof
	As in the proof of Theorem \ref{th-Rate}, denote $F_k = F(x_k) - F^*$ and $g_k = \| F'(x_k) \|_*$. 
	Then, we have                                                                                                                                                                                                                                                                                                                                                                                                                                                                                                                                                                                                                                                                                                                                                                                                                                                                                                                                                                                                                                                                                                                                                                                                                                                                                                                                                          
	$$
	\ba{rcl}
	\ln {1 \over F_{k+1}} - \ln {1 \over F_k} & = & \ln \left( 1 + {F_k - F_{k+1} \over F_{k+1}} \right) \; \refGE{eq-Prog} \; \ln \left( 1 + {\sqrt{3} g_{k+1}^2 \over 2 \sqrt{H} c^2 g_k^{1/2} F_{k+1}} \right)\\
	\\
	& \refGE{eq-UBound} & \ln \left( 1 +  {3 \over 4c^2} \sqrt{3 \sigma_3 \over H} \cdot {g_{k+1}^{1/2} \over g_k^{1/2}} \right) \; = \; \ln \left( 1 + S \cdot \sqrt{g_{k+1} \over c g_k } \right),
	\ea
	$$
	where $S = {3 \over 4c^{3/2}} \sqrt{3 \sigma_3 \over H}$. Denote $\tau_k = \sqrt{g_{k+1} \over c g_k } \refLE{eq-GNorm} 1$. Since $\ln(\cdot)$ is a concave function, we have $\ln(1+S \tau_k) \geq \tau_k \ln(1+S)$. Hence,
	$$
	\ba{rcl}
	\xi_k \Def \ln {g_0 D \over F_k} \; \geq \; \ln {F_0 \over F_k} & \geq & \ln(1+S) \sum\limits_{i=0}^{k-1} \tau_i \; \geq \;
	{k \over c^{1/2}} \ln(1+S) \left(\prod\limits_{i=0}^{k-1} {g_{i+1}^{1/2} \over g_i^{1/2}}\right)^{1/k} \\
	\\
	&  = & {k \over c^{1/2}} \ln(1+S) \left({g_k \over g_0}\right)^{1/(2k)}.
	\ea
	$$
	Note that $\left({g_k \over g_0}\right)^{1/(2k)} = \exp\left(- {1 \over 2k} \ln {g_0 \over g_k} \right) \geq 1 + {1 \over 2k} \ln {g_k \over g_0} \geq 1 + {1 \over 2k} \ln {F_k \over g_0 D}  = 1 - {1 \over 2k} \xi_k$. Thus,
	$$
	\ba{rcl}
	\xi_k & \geq & {k \ln(1+S) \over c^{1/2} + \half \ln(1+S)},
	\ea
	$$
	and this is inequality (\ref{eq-RateUC3}).
	\qed
	\BR
	in accordance to the estimate (\ref{eq-RateUC3}), the highest rate of convergence corresponds to the maximal value of $S$. This means that we need to minimize the factor $c^{3/2} H^{1/2}$ in $H$. The optimal value is given by $H_{\#} = {3 \sigma}L_2$. In this case,
	\beq\label{eq-LRate}
	\ba{rcl}
	S & = & {\sigma } \sqrt{\sigma_3 \over 6 L_2} \; > \; 0.4 \sigma \sqrt{\sigma_3 \over L_2}.
	\ea
	\eeq
	\ER
	
	Finally, let us prove a superlinear rate of local convergence for the scheme (\ref{met-GNewt}).
	\BT\label{th-Super}
	Let function $f(\cdot)$ be strongly convex on $\dom \psi$ with parameter $\mu > 0$. If $H \geq L_2$, then, for any $k \geq 0$ we have
	\beq\label{eq-LocRate}
	\ba{rcl}
	\| F'(x_{k+1} ) \|_* & \leq & {2 c\over \mu} \sqrt{H \over 3} \| F'(x_k) \|_*^{3/2}.
	\ea
	\eeq
	\ET
	\proof
	Indeed, for any $k \geq 0$ we have
	$$
	\ba{rcl}
	{\mu \over 2} \| x_{k+1} - x_k \|^2 & \leq & \half \la \nabla^2f(x_k)(x_{k+1} - x_k), x_{k+1} - x_k \ra \\
	\\
	& \refLE{eq-DecF} & F(x_k) - F(x_{k+1}) \; \leq \; \| F'(x_k) \|_* \| x_k - x_{k+1} \|.
	\ea
	$$
	Therefore,
	$$
	\ba{rcl}
	\| F'(x_{k+1}) \|_* & \refLE{eq-GNorm} & \sigma c A_k \| x_{k+1} - x_k \| \; \leq \; {2 \sigma c \over \mu} A_k \| F'(x_k) \|_*\\
	\\
	& \refEQ{def-AH} & {2 c \over \mu} \sqrt{H \over 3} \| F'(x_k) \|_*^{3/2}. \QR
	\ea
	$$
	
	Thus, the region of superlinear convergence of method (\ref{met-GNewt}) is as follows:
	\beq\label{eq-Region}
	\ba{c}
	{\cal R}_Q \Def \left\{ x \in \dom \psi:\; \| F'(x) \|_* \leq {3 \mu^2 \over 4 H c^2} \right\}.
	\ea
	\eeq
	Note that outside this region, the constant of strong convexity of the objective function in problem (\ref{def-TA}) with $A = A_H(x)$ satisfies the following lower bound:
	\beq\label{eq-MuLow}
	\ba{rcl}
	\sigma A_H(x) & \geq & {\mu \over 2c}, \quad x \not\in {\cal R}_Q.
	\ea
	\eeq

	\section{Estimating the norm of the gradient}\label{sc-Grad}
	\SetEQ
	
	Let us estimate the efficiency of method (\ref{met-GNewt}) in decreasing the norm of gradients. For that, we are going to derive an upper bound for the number of steps $N$ of method (\ref{met-GNewt}), for which we still have
	\beq\label{eq-GGap}
	\| F'(x_k) \|_* \geq \delta > 0, \quad 0 \leq k \leq N.
	\eeq
	In this section, we use notation of Section \ref{sc-MinProc}:
	$$
	\ba{rcl}
	F_k & = & F(x_k) - F^*, \quad g_k \; = \; \| F'(x_k) \|_*.
	\ea
	$$
	
	Firstly, consider the case when the smooth component $f(\cdot)$ in the objective function of problem (\ref{prob-Uni}) satisfies condition (\ref{eq-Lip}). Then
	\beq\label{eq-Prog1}
	\ba{rcl}
	F_k - F_{k+1} & \refGE{eq-Prog} & \kappa {g_{k+1}^2 \over g_k^{1/2}}, \quad \kappa \Def {1 \over 2c^2} \sqrt{3 \over H}.
	\ea
	\eeq
	It is convenient to assume that the number of iteration $N$ of the method is a multiple of three:
	\beq\label{eq-Even}
	N = 3m, \quad m \geq 1.
	\eeq
	Then for the last $m$ iterations of the scheme we have
	\beq\label{eq-Stage2}
	\ba{rcl}
	F_{2m} & \geq & F_{2m} - F_{3m} \; \geq \; \kappa \sum\limits_{i=0}^{m-1} {g_{2m+i+1}^2 \over g_{2m+i}^{1/2}} \; \refGE{eq-GGap} \; \kappa \delta^{3/2} \sum\limits_{i=0}^{m-1} {g_{2m+i+1}^{1/2} \over g_{2m+i}^{1/2}}\\
	\\
	& \geq & \kappa m \delta^{3/2} \left({g_{3m}^{1/2} \over g_{2m}^{1/2}}\right)^{1/m} \; \refGE{eq-GGap} \; \kappa m \delta^{3/2} \left({\delta^{1/2} \over g_{2m}^{1/2}}\right)^{1/m}.
	\ea
	\eeq
	At the same time, for the first $2m$ iterations we obtain
	\beq\label{eq-Stage1}
	\ba{rcl}
	{1 \over F_{2m}^{1/2}} - {1 \over F_0^{1/2}} & \refGE{eq-Int} & {2 m \over 4 c^2} \sqrt{3 \over H D^3} \left({F_{2m}g_{2m}^{1/2} \over F_0 g_0^{1/2}} \right)^{1/(2m)} \; = \; \kappa m D^{-3/2}\left({F_{2m}g_{2m}^{1/2} \over F_0 g_0^{1/2}} \right)^{1/(2m)}.
	\ea
	\eeq
	Hence, using inequality (\ref{eq-Stage2}) and squared inequality (\ref{eq-Stage1}), we obtain the following:
	$$
	\ba{rcl}
	1 & \geq & \left(1 - \sqrt{F_{2m} \over F_0} \right)^2 \; = \; \left( {1 \over F_{2m}^{1/2}} - {1 \over F_0^{1/2}} \right)^2 \cdot F_{2m} \; \geq \;
	\left({ \kappa m \delta^{1/2} \over D} \right)^3 \left({F_{2m} \delta^{1/2} \over F_0 g_0^{1/2}} \right)^{1/m}
	\ea
	$$
	Note that $g_{2m} \refLE{eq-GNorm} c^{2m} g_0$. Therefore, 
	$$
	\ba{rcl}
	F_{2m} & \refGE{eq-Stage2} & \kappa m \delta^{3/2} \left({\delta^{1/2} \over c^m g_0^{1/2}}\right)^{1/m},
	\ea
	$$
	and we obtain
	$$
	\ba{rcl}
	1 & \geq & \left({ \kappa m \delta^{1/2} \over D} \right)^3 \left({ \kappa m \delta^{2} \over c F_0 g_0^{1/2}} \cdot \left({\delta^{1/2} \over g_0^{1/2}}\right)^{1/m}\right)^{1/m} \\
	\\
	& \geq & \left({\kappa m \delta^{1/2} \over D} \right)^{3 + {1 \over m}} \left({\delta^{1/2} \over g_0^{1/2}}\right)^{ (3 + {1 \over m}) {1 \over m}}
	\Big( c \Big)^{-{1 \over m}}. 
	\ea
	$$
	Thus, we can prove the following theorem.
	\BT\label{th-CompG}
	Under condition (\ref{eq-GGap}), the number of steps of method (\ref{met-GNewt}) satisfies the following bound:
	\beq\label{eq-CompG}
	\ba{rcl}
	N & \leq & 2 c^2 \sqrt{3 H D^2 \over \delta} + {3 \over 2} \ln {g_0 \over \delta} + \ln c.
	\ea
	\eeq
	\ET
	\proof
	Indeed,
	$$
	\ba{rcl}
	1 & \geq & {\kappa m \delta^{1/2} \over D} \Big( {\delta \over g_0} \Big)^{1 \over 2m} \Big(  c \Big)^{-{1 \over 3m+1}} \; = \;  {\kappa m \delta^{1/2} \over D} \exp \left( - {1 \over 2m} \ln \left[ {g_0 \over \delta} \Big( c \Big)^{{2 m \over 3 m+1}} \right] \right)\\
	\\
	& \geq & {\kappa  \delta^{1/2} \over D} \left( m - \half \ln {g_0 \over \delta} - {m \over 3m+1} \ln c \right) \; \geq \; {\kappa  \delta^{1/2} \over D} \left( m - \half \ln {g_0 \over \delta} - {1 \over 3} \ln c \right),
	\ea
	$$
	and this is inequality (\ref{eq-CompG}).
	\qed
	
	Finally, let us estimate the efficiency of method (\ref{met-GNewt}) under additional assumption of uniform convexity (\ref{def-UC3}). From the proof of Theorem \ref{th-UC3}, we know that
	$$
	\ba{rcl}
	\ln{F_0 \over F_{2m}} & \geq & {2m \over c^{1/2}} \ln (1+S) \left( {g_{2m} \over g_0} \right)^{1/(2m)} \; \geq \; {2 m \over c^{1/2}} \ln (1+S) \exp \left( - {1 \over 2m} \ln {g_{0}\over g_{2m}} \right)\\
	\\
	& \geq & {1 \over c^{1/2}} \ln (1+S) \left( 2m  - \ln {g_{0}\over g_{2m}} \right) \; \refGE{eq-GGap} \; {1 \over c^{1/2}} \ln (1+S) \left( 2m  - \ln {g_{0}\over \delta} \right).
	\ea
	$$
	On the other hand,
	$$
	\ba{rcl}
	\ln F_{2m} & \refGE{eq-Stage2} &  \ln (\kappa m \delta^{3/2}) + {1 \over 2m} \ln {\delta \over g_{2m}} \; \refGE{eq-GNorm} \;  \ln (\kappa m \delta^{3/2}) + {1 \over 2m} \ln {\delta \over g_0} - \ln c .
	\ea
	$$
	Thus,
	$$
	\ba{rcl}
	\ln (cF_0) & \geq &  {2m \over c^{1/2}} \ln (1+S) - {1 \over c^{1/2}} \ln (1+S) \ln {g_{0}\over \delta} + \ln (\kappa m \delta^{3/2}) + {1 \over 2m} \ln {\delta \over g_0}.
	\ea
	$$
	In other words,
	$$
	\ba{rcl}
	\ln {c F_0 \over \kappa g_0^{3/2}} & \geq &  {2m \over c^{1/2}} \ln (1+S)
	- {1 \over c^{1/2}} \ln (1+S) \ln {g_{0}\over \delta} + {3 \over 2} \ln {\delta \over g_0}  - \ln {1 \over m} + {1 \over 2m} \ln {\delta \over g_0}\\
	\\
	& = & {2m \over c^{1/2}} \ln (1+S)
	-  \left[ {1 \over 2m} +  {1 \over c^{1/2}} \ln (1+S)  + {3 \over 2} \right] \ln {g_{0}\over \delta}  - \ln {1 \over m}.
	\ea
	$$
	Thus, we have proved the following theorem.
	\BT\label{th-CompGUC}
	Under condition (\ref{eq-GGap}) and uniform convexity (\ref{def-UC3}), the number of steps of method (\ref{met-GNewt}) satisfies the following bound:
	\beq\label{eq-CompG2}
	\ba{rcl}
	N & \leq & {3 c^{1/2} \over \ln(1+S)} \left\{  \ln {c F_0 \over \kappa g_0^{3/2}} + \left[ {1 \over 2m} +  {1 \over c^{1/2}} \ln (1+S)  + {3 \over 2} \right] \ln {g_{0}\over \delta}  \right\}\\
	\\
& \refLE{eq-UBound} & {3 c^{1/2} \over \ln(1+S)} \ln { 3 c F_0 \over 2 \kappa 
\sqrt{\sigma_3}}
+ 3 \Big [ 1 + {2 c^{1/2} \over \ln(1+S)} \Big] \ln {g_{0}\over \delta}.
	\ea
	\eeq
	\ET

	\section{Adaptive line search}\label{sc-LS}
	\SetEQ
	
	The main advantage of the method ({\ref{met-GNewt}) consists in its easy implementation. Indeed, in the case $\psi(\cdot) \equiv 0$ with $\dom \psi = \E$, the iteration (\ref{def-TA}) is reduced mainly to matrix inversion, the very standard operation of Linear Algebra, which is available in the majority of software packages. However, for the better performance of this scheme, it is necessary to apply a dynamic strategy for updating the step-size coefficient $H$. Let us show how this can be done.
		
		Consider the following optimization method.
		\beq\label{met-GNewtF}
		\ba{|c|}
		\hline \\
		\mbox{\bf Gradient Regularization of Newton Method with Line Search}\\
		\\
		\hline \\
		\ba{l}
		\mbox{{\bf Initialization.} Choose $H_0 \leq L_2$, $x_0 \in \dom \psi$, and $F'_0 \in \partial F(x_0)$.}\\
		\\
		\mbox{{\bf $k$th iteration ($k \geq 0$).} 1). Set $g_k = \| F'_k \|_*$.}\\
		\\
		\mbox{2). Find the least $i = i_k \geq 0$: s.t. for $H = 2^i H_k$ and $T = T_{A_H(x_k)}(x_k)$, we}\\
		\\
		\mbox{have $f(T) \leq f(x_k) + \la \nabla f(x_k), T - x_k \ra + \half \nabla^2f(x_k)[T-x_k]^2 + {H \over 6} \| T - x_k \|^3$.}\\
		\\
		\mbox{3). Set $A_k = {1 \over \sigma}\sqrt{{2^{i_k} \over 3} H_k g_k}$, $x_{k+1} = T_{A_k}(x_k)$, $H_{k+1} = \max\{ H_0, 2^{i_k-1} H_k\}$, and}\\
		\\
		F'_{k+1} = \nabla f(x_{k+1}) - \nabla f(x_k) - \nabla^2f(x_k)(x_{k+1}-x_k) - A_k(\nabla d(x_{k+1}) - \nabla d(x_k))\\
		\ea\\
		\\
		\hline
		\ea
		\eeq
		
	Note that this scheme does not depend on any particular value of the Lipschitz constant. By definitions of the updates and from inequality~\eqref{eq-FLip}, we conclude that inequalities 
$H_0 \leq H_k \leq L_2$ and $2^{i_k} H_k \leq 2L_2$ imply $H_{k+1} \leq L_2$. Thus, 
	\beq \label{eq-HBound}
	\ba{rcl}
	H_0 & \leq & H_k \;\; \leq \;\; L_2, \quad 2^{i_k} H_k \; \leq \; 2L_2, \qquad k \geq 0. 
	\ea
	\eeq
	Hence, from Theorem~\ref{th-Main}, we have the following progress
	established for each iteration $k \geq 0$:
	$$
	\ba{rcl}
	F(x_k) - F(x_{k + 1}) & \geq & \frac{1}{2c_0^2} \sqrt{\frac{3}{2L_2}} 
	\cdot \frac{\| F'(x_{k + 1}\|_*^2}{\|F'(x_k)\|_*^{1/2}},
	\ea
	$$ 
	where
	$$
	\ba{rcl}
	c_0 & \Def &  \sigma^{-1} + \frac{3L_2}{2H_0} .
	\ea
	$$
	Repeating the reasoning of Theorem~\ref{th-Rate}, we obtain the following complexity result.
	\BT \label{th-LS}
	Let $F(x_k) - F^* \geq \epsilon$ for some iteration $k \geq 0$ of method~\eqref{met-GNewtF}. Then, 
	$$
	\ba{rcl}
		k & \leq & 4c_0^2 \sqrt{ \frac{2L_2 D^3}{3\epsilon} }
		+ \ln \frac{(F(x_0) - F^{*})\|F'(x_0)\|_*^{1/2} D^{1/2}}{\epsilon^{3/2}}. \QF
	\ea
	$$
	\ET

	\newpage
	\section{Acceleration}\label{sc-Acc}
	\SetEQ

	Let us present a conceptual acceleration of our method,
	that is based on the contracting proximal iterations~\cite{doikov2020contracting}.
	
	First, we fix an auxiliary prox-function $\phi(\cdot)$ that we assume to be uniformly 
	convex of degree three with respect to the initial norm:
	\beq \label{eq-AccProx}
	\ba{rcl}
	\beta_{\phi}(x, y) & = & \phi(y) - \phi(x) - \la \nabla \phi(x), y - x \ra
	\;\; \geq \;\; \frac{1}{3}\|y - x\|^3, \qquad \forall x, y \in \dom \psi.
	\ea
	\eeq
	At each iteration $k \geq 0$ of the accelerated scheme, we form the following functions:
	$$
	\ba{rcl}
	g_{k + 1}(x) & \Def & B_{k + 1} f\Bigl( \frac{b_{k + 1} x + B_k x_k}{B_{k + 1}} \Bigr), \\
	\\
	h_{k + 1}(x) & \Def & g_{k + 1}(x) + b_{k + 1} \psi(x) + \beta_{\phi}(v_k; x),
	\ea
	$$
	where $\{ b_k \}_{k \geq 1}$ is a sequence of positive numbers, $B_k \Def \sum\limits_{i = 1}^k b_i$, $B_0 \Def 0$, and
	$$
	\ba{c}
	\{ x_k \}_{k \geq 0}, \quad \{ v_k \}_{k \geq 0}, \quad x_0 = v_0,
	\ea
	$$
	are sequences of trial points that belong to $\dom \psi$.
	
	Note that the derivatives of $g_{k + 1}( \cdot )$ and $f(\cdot)$ are related as follows:
	$$
	\ba{rcl}
		D^3 g_{k + 1}(x) & \equiv & 
		\frac{b_{k + 1}^3}{B_{k + 1}^2} D^3 f\Bigl( \frac{b_{k + 1} x + B_k x_k}{B_{k + 1}} \Bigr).
	\ea
	$$
	For simplicity of the presentation, we assume that $f$ is three times differentiable
	on the open set containing $\dom \psi$.
	Let us choose
	$$
	\ba{rcl}
	b_{k} & := & \frac{k^2}{9 L_2(f)}.
	\ea
	$$
	Then, $B_k = \frac{1}{9 L_2(f)} \sum\limits_{i = 1}^k i^2
	\geq \frac{k^3}{27 L_2(f)}$. Therefore, for any $h \in \E$:
	$$
	\ba{rcl}
	|D^3 g_{k + 1}(x)[h]^3| & \leq & \frac{1}{L_2(f)}|
	D^3 f\Bigl( \frac{b_{k + 1} x + B_k x_k}{B_{k + 1}} \Bigr)|
	\;\; \leq \;\; \|h\|^3,
	\ea
	$$
	thus $L_2(g_{k + 1}) = 1$, and we can minimize objective $h_{k + 1}$
	very efficiently by using our method~\eqref{met-GNewt}. Namely,
	in order to find a point $v$ with a small norm of a subgradient:
	$$
	\ba{rcl}
	\| g \|_{*} & \leq & \delta, \qquad g \in \partial h_{k + 1}(v),
	\ea
	$$
	the method needs to do no more than
	$$
	\ba{rcl}
	N & \overset{\eqref{eq-CompG}}{\leq} & \tilde{O}( \ln \frac{1}{\delta} ) 
	\ea
	$$
	steps, where $\tilde{O}(\cdot)$ hides absolute constants and logarithmic factors
	that depends on the initial  residual and subgradient norm.
	
	Let us write down the accelerated algorithm.
		\beq\label{met-Accel}
	\ba{|c|}
	\hline \\
	\mbox{\bf Acceleration of Newton Method with Gradient Regularization}\\
	\\
	\hline \\
	\ba{l}
	\mbox{{\bf Initialization.} Choose $x_0 \in \dom \psi$ and $\delta > 0$. 
		Set $v_0 = x_0$, $B_0 = 0$.}\\
	\\
	\mbox{{\bf $k$th iteration ($k \geq 0$).} 1). Set $b_{k + 1} = \frac{k^2}{9L_2(f)}$
		and $B_{k + 1} = B_k + b_{k + 1}$.}\\
	\\
	\mbox{2). Form the auxiliary objective $h_{k + 1}(\cdot)$. Find a point $v_{k + 1}$ by method~\eqref{met-GNewt} } \\
	\\
	\mbox{such that $\; \|g\|_{*} \; \leq \; \delta \;$ for some 
		$\;g \in \partial h_{k + 1}(v_{k + 1})$.} \\
	\\
	\mbox{3). Set $x_{k + 1} = \frac{b_{k + 1} v_{k + 1} + B_k x_k}{B_{k + 1}}$.} \\
	\ea\\
	\\
	\hline
	\ea
	\eeq
	
	Applying directly Theorem 3.2 and the corresponding Corollary 3.3 
	from~\cite{doikov2020contracting}, we get the following complexity bound.
	
	\BT Let us set $\delta = \frac{1}{2 \cdot 3^{7/3}} \cdot \bigl(\frac{\epsilon}{L_2(f)}\bigr)^{2/3}$
	in method~\eqref{met-Accel}, and let
	$$
	\ba{rcl}
	k & = & \Bigl\lceil \bigl( 2 \cdot 3^3 \bigr)^{1/2} \cdot 
	\Bigl( \frac{L_2(f) \beta_\psi(x_0; x^{*})}{\epsilon} \Bigr)^{1/3}    \Bigr\rceil.
	\ea
	$$
	Then, $F(x_k) - F^{*} \leq \epsilon.  \QR$
	\ET

	\bibliographystyle{plain}

\end{document}